\documentclass[10pt,twoside]{article}

\usepackage{graphicx}

\usepackage{amssymb}
\usepackage{latexsym}

\makeatletter

\@addtoreset{equation}{section} \makeatother
\newtheorem{theorem}{Theorem}[section]
\newtheorem{remark}[theorem]{Remark}
\newtheorem{lemma}[theorem]{Lemma}
\newtheorem{proposition}[theorem]{Proposition}
\newtheorem{corollary}[theorem]{Corollary}
\newtheorem{definition}[theorem]{Definition}
\newtheorem{example}[theorem]{Example}

\def\1{\mathbf{1}}
\def\:{\lrcorner}
\def\#{\sharp}

\def\inv#1{\raise.1em\hbox to 0pt{$^{-1}$\hss}_{#1}\;}

\newcommand{\bean}{\begin{eqnarray*}}
\newcommand{\eean}{\end{eqnarray*}}
\newcommand{\benu}{\begin{enumerate}}
\newcommand{\eenu}{\end{enumerate}}
\newcommand{\eea}{\end{eqnarray}}
\newcommand{\bea}{\begin{eqnarray}}

\newcommand{\q}{p'}     

\newcommand{\be}{\begin{equation}}
\newcommand{\ee}{\end{equation}}

\newcommand{\N}{{\mathbb N}}

\newcommand{\R}{{\mathbb R}}
\newcommand{\LL}{{\mathbb L}}

\newcommand{\noi}{\noindent}
\newcommand{\sm}{\smallskip}

\newcommand{\ben}{\begin{enumerate}}
\newcommand{\een}{\end{enumerate}}
\newcommand{\bit}{\begin{itemize}}
\newcommand{\eit}{\end{itemize}}
\newcommand{\edoc}{\end{document}}

\newcommand{\bdefi}{\begin{definition}}
\newcommand{\btheo}{\begin{theorem}}
\newcommand{\bprop}{\begin{proposition}}
\newcommand{\brema}{\begin{remark}}
\newcommand{\bcoro}{\begin{corollary}}
\newcommand{\blemm}{\begin{lemma}}
\newcommand{\bexam}{\begin{example}}

\newcommand{\edefi}{\end{definition}}
\newcommand{\etheo}{\end{theorem}}
\newcommand{\eprop}{\end{proposition}}
\newcommand{\erema}{\end{remark}}
\newcommand{\ecoro}{\end{corollary}}
\newcommand{\elemm}{\end{lemma}}
\newcommand{\eexam}{\end{example}}

\newcommand{\proof}{\noindent {\em Proof.\ }}
\newcommand{\cvd}{\ \rule{0.5em}{0.5em} \sm \noi}

\title{ {\bf Lorentzian manifolds \\ isometrically embeddable in  $\LL^N$
 \vspace{0.5in}  }}

\author{\bf O. M\"uller$^{1}$, M.
S\'anchez$^{2}$
\vspace{1mm}\\
$^1$ {\it\small Instituto de Matem\'aticas,}\\
{\it \small UNAM Campus Morelia,}\\
{\it \small C. P. 58190 Morelia, Michoac\'an, M\'exico.}\\
{\it \small email: olaf@matmor.unam.mx}
\vspace{1mm} \\
$^2$ {\it\small Departamento de Geometr\'{\i}a y Topolog\'{\i}a}\\
{\it\small Facultad de Ciencias, Universidad de Granada}\\
{\it\small Campus Fuentenueva s/n, 18071 Granada, Spain}}
\begin{document}

\parindent=5mm
\date{}
\maketitle

\begin{quote}

\noindent {\small \bf Abstract.} {\small In this article,  the
Lorentzian manifolds isometrically embeddable in $\LL^N$ (for some
large $N$, in the spirit of Nash's theorem) are characterized as a
a subclass of the set of all stably  causal spacetimes;
concretely, those which admit a smooth time function $\tau$ with
$|\nabla \tau|>1$. Then, we prove that any globally hyperbolic
spacetime $(M,g)$ admits such a function, and, even more, a global
orthogonal decomposition $M=\R \times S, g=-\beta dt^2 + g_t$ with
bounded function $\beta$ and Cauchy slices.

In particular,  a proof of a  result stated by C.J.S. Clarke is
obtained: any globally hyperbolic spacetime
 can be isometrically embedded in Minkowski spacetime
$\LL^N$. 
The role of the so-called ``folk problems on smoothability'' in
Clarke's approach is also discussed.
}\\

\end{quote}

\begin{quote}

\end{quote}
\begin{quote}


\end{quote}

\newpage



\section{Introduction}

A celebrated  theorem by J. Nash \cite{Na} states that any $C^3$
Riemannian manifold can be isometrically embedded in any open
subset of some Euclidean space $\R^{N}$ for large $N$.  Greene
\cite{Gr} and Clarke \cite{Cl} showed independently, by means of
simple algebraic reasonings, that Nash's theorem can be extended
to indefinite (even degenerate) metrics, that is, any
semi-Riemannian manifold can be smoothly isometrically embedded in
any open subset of semi-Euclidean space $\R^N_s$ for large enough
dimension $N$ and index $s$. Moreover, they also reduced the Nash
value for $N$:
 Greene by using
 the implicit function theorem by Schwartz \cite{Sc}, and
Clarke by means of a  technique inspired in Kuiper's \cite{Ku},
which yields  $C^k$ isometric embeddings with $3\leq k <\infty$.

Nevertheless, a new problem appears when a semi-Riemannian
manifold of index $s$ is going to be embedded in a semi-Euclidean
space of the same index $\R^N_s$. We will focus on the simplest
case $s=1$, i.e., the isometric embedding of a Lorentzian manifold
$(M,g)$ in Minkowski spacetime $\LL^N$. Such an embedding will not
exist in general;  for example, the existence of a causal closed
curve in $M$ contradicts the possibility of an embedding in
$\LL^N$. So, the first task is to characterize the class of
isometrically embeddable spacetimes. This is the role of our first
result (Section \ref{s3}):

\btheo \label{t1} Let $(M,g)$ be a  Lorentzian manifold. The
following assertions  are equivalent:

(i) $(M,g)$ admits a  isometric embedding in $\LL^N$ for some
$N\in \N$.

(ii) $(M,g)$ is a stably causal spacetime with a {\em steep
temporal function}, i.e., a smooth function $\tau$ such that
$g(\nabla \tau, \nabla \tau)\leq -1$.

\etheo

Again, this theorem is carried out by using some simple arguments,
which essentially reduce the hardest problem to the Riemannian
case. So, this result (and the subsequent ones on isometric
embeddings) is obtained under the natural technical conditions
which comes from the Riemannian setting: (a) $(M,g)$ must be $C^k$
with $3 \leq k \leq \infty$, and all the other elements will be as
regular as permitted by $k$, and (b) the smallest value of $N$  is
$N=N_0(n)+1$, where $n$ is the dimension of $M$ and $N_0(n)$ is
the optimal bound in the Riemannian case (see \cite{HH} for a
recent summary on this bound). We will not care about the local
problem (see \cite{Gr}, a summary in Lorentzian signature can be
found in \cite{S}); recall also that, locally, any spacetime
fulfills condition (ii). So, the main problem we will consider
below, is the existence of a global steep temporal function as
stated in (ii).

It is known that any stably causal spacetime admits a {\em time}
function, which can be smoothed into a {\em temporal} one $\tau$
(see Section \ref{s2} for definitions and background).
Nevertheless, the condition of being {\em steep}, $|\nabla
\tau|\geq 1$ cannot be fulfilled for all stably causal spacetimes.
In fact, a simple counterexample, which works even in the causally
simple case, is provided below (Example \ref{ex}). Notice that
causal simplicity is the level in the standard {\em causal
hierarchy of spacetimes} immediately below global hyperbolicity.
So, the natural question is to wonder if any globally hyperbolic
spacetime admits a steep temporal function $\tau$.

The existence of  embeddings in $\LL^N$ for globally hyperbolic
spacetimes  was stated by Clarke \cite[Sect. 2]{Cl}. In his
approach, a function $f: M \rightarrow \LL^2$ with a similar role
to the steep temporal function above is used. Nevertheless, as in
other papers of that epoch, his construction of $f$ is affected by
the so-called ``folk problems'' of smoothability of
causally-constructed functions. So, as will be discussed in the
Appendix, if Clarke's proof is completed, then a new type of
causally-constructed functions will be shown to be smooth (or at
least smoothable).

Apart from the consequence of the  embedding in $\LL^N$, the
existence of  a steep temporal $\tau$ is relevant for the
structure of globally hyperbolic spacetimes. In fact, both
questions, problems of smoothability and structure of globally
hyperbolic spacetimes, are linked since Geroch's landmark about
topological splittings \cite{Ge}. More precisely, recently
 any globally
hyperbolic spacetime $(M,g)$ has been proved to admit a {\em
Cauchy orthogonal decomposition} \be\label{cosplit} M=\R\times S ,
\quad \quad g=-\beta d{\cal T}^2 + g_{\cal T} , \ee where
$\beta>0$ is a function on $M$, $g_{\cal T}$ is a Riemannian
metric on $S_{\cal T}:=\{{\cal T}\}\times S$ smoothly varying with
${\cal T}$, and each slice $S_{\cal T}$
 becomes a Cauchy hypersurface \cite{BeSa2}.
This result, which improves  Geroch's topological splitting $M
\cong \R\times S$ as both, a differentiable and orthogonal one, is
proved by showing that, starting at Geroch's Cauchy time function,
one could obtain a Cauchy  smooth time function with timelike
gradient. Now, recall that, if this {\em Cauchy temporal} function
is steep, then the function $\beta$ is upper bounded by one (Lemma
\ref{l1}).  One of us suggested possible analytical advantages of
a strengthened decomposition (\ref{cosplit}),  where additional
conditions on the elements $\beta, g_{\cal T}$ are imposed
\cite{Mu}. In particular, such a decomposition is called there a
{\em b-decomposition} if the function $\beta$ (the {\em lapse} in
relativist's terminology) is bounded. Our next result is then
(Section \ref{s4}):

\btheo \label{t2} Any globally hyperbolic spacetime admits a steep
Cauchy temporal function ${\cal T}$ and, so, a Cauchy orthogonal
decomposition (\ref{cosplit}) with (upper) bounded  function
$\beta$.
 \etheo

\brema \label{r}{\em From the technical viewpoint, the
decomposition (\ref{cosplit}) was carried out in \cite{BeSa2} by
proving the existence of a Cauchy temporal function; moreover, a
simplified argument shows the existence of a temporal function in
any stably causal spacetime (\cite{BeSa2}, see also the discussion
in \cite{Sa-bras}). Our proof is completely self-contained, as it
re-proves the existence of the Cauchy temporal function with
different and somewhat simpler arguments, as well as a stronger
conclusion. Nevertheless, we use some technical elements
(remarkably, Proposition \ref{lemilla}) which hold in the globally
hyperbolic case, but not in the stably causal one\footnote{Notice
also that only $C^1$ differentiability is needed for these
results.}.}\erema

Summing up, we   emphasize the following consequences of
previous theorems (for the second one recall also Proposition
\ref{p0}).

\bcoro \label{c1} (1)  Any  globally hyperbolic spacetime can be
isometrically embedded in some $\LL^N$.

(2) A Lorentzian manifold is a stably causal spacetime if and only
if it admits a conformal embedding in some $\LL^N$. In this case,
there is a representative of its conformal class whose
time-separation (Lorentzian distance) function is finite-valued.

\ecoro  Notice also that, as an immediate consequence, a stably
causal spacetime is not globally hyperbolic if and only if it is
conformal to a spacetime non-isometrically embeddable in $\LL^N$
(see Example \ref{ex}).

After some preliminaries in the next section, Sections \ref{s3},
\ref{s4} are devoted, respectively, to prove  Theorems \ref{t1}
and \ref{t2}, as well as to discuss its optimality and
consequences. From the technical viewpoint, it is worth pointing
out the introduction of two elements in the first part of Section
\ref{s4}: a semi-local temporal function for subsets of type
$J^\pm(p)\cap J^\mp(S)$ (Prop. \ref{lemilla}) and {\em fat cone
coverings} for any Cauchy hypersurface $S$ (Prop. \ref{fatcone}).
Finally, in the Appendix, Clarke's technique for globally
hyperbolic spacetimes is discussed, and new causal problems on
smoothability, which may have their own interest, are suggested.

\section{Preliminaries}\label{s2}

In what follows, any  semi-Riemannian manifold will be  $C^k$,
with $3 \leq k \leq \infty$ as in Nash's theorem, and  will be
assumed to be connected without loss of generality. Any geometric
element on the manifold will be called smooth if it has the
highest order of differentiability allowed by $k$. For an
immersion $i: M \rightarrow \bar M$ only injectivity of each
$di_p, p\in M$ is required;  the injectivity of  $i$, as well as
being a homeomorphism onto its image, are required additionally
for $i$ to be an embedding.

Our notation and conventions on causality will be standard as, for
example, in \cite{BEE} or \cite{O}. Nevertheless, some terminology
on the solution of the so-called ``folk problems of
smoothability'' introduced in \cite{BeSa1, BeSa2} are also used
here (see \cite{MS} for a  review). In particular, a Lorentzian
manifold $(M,g)$ is a  manifold $M$ endowed with a metric tensor
$g$ of index one
 $(-,+,\dots,+)$,
a tangent vector $v\in T_pM$ in $p\in M$, is timelike (resp.
spacelike; lightlike; causal)  when $g(v,v)<0$ (resp, $g(v,v)>0$;
 $g(v,v)=0$ but $v\not=0$; $v$ is timelike or lightlike); so, following \cite{MS},
the vector 0 will be regarded as non-spacelike and non-causal --
even though this is not by any means the unique convention in the
literature. For any vector $v$, we write $\vert v \vert := \sqrt{
\vert g(v, v) \vert }$. A spacetime is a time-orientable
Lorentzian manifold, which will be assumed to be time-oriented
(choosing any of its two time-orientations) when necessary; of
course, the choice of the time-orientation for submanifolds
conformally immersed in $\LL^N$ will agree with the induced from
the canonical time-orientation of $\LL^N$. The associated {\em
time-separation} or {\em Lorentzian distance} function will be
denoted by $d$, $d(p,q) := \rm{sup}_{c \in \Omega(p,q)} l(c) $
where the supremum is taken over the space $\Omega(p,q)$ of
future-directed causal $C^1$ curves from $p$ to $q$ parametrized
over the unit interval (if this space is empty, $d$ is defined
equal to $0$), and $ l(c) := \int_0^1 \vert \dot{c} (t) \vert dt$
for such a curve. The following elements of causality must be
taken into account (they are explained in detail in \cite{MS}).

\bit \item A {\em time function} $t$ on a spacetime is a
continuous function which increases strictly on any
future-directed causal curve.  Recently \cite{BeSa2}, it has been
proved that this is also equivalent to the existence of a {\em
temporal function} $\tau$, i.e., a smooth time function with
everywhere past-directed timelike gradient $\nabla \tau$. This
also ensures the  folk claim  that, for a spacetime, the existence
of  a time funtion is equivalent to be {\em stably causal} (i.e.,
 if the lightcones of the spacetime  are slightly opened then it
remains causal), see  \cite[Fig. 2, Th. 4.15, Rem. 4.16]{Sa-bras}
or \cite[Th. 3.56]{MS}. Along the present paper, a temporal
function will be called {\em steep} if $|\nabla \tau|\geq 1$; as
we will see, not all stably causal spacetimes admit a steep
temporal function.

\item After stable causality, the two next steps in the so-called
{\em causal ladder} or {\em causal hierarchy} of spacetimes are:
causal continuity (the volume functions $t^\pm(p)= \mu(I^\pm(p)),
p\in M$ are time functions for one, and then for all,
 measure associated to any auxiliary semi-Riemannian metric such
that $\mu(M)<\infty$) and causal simplicity (the spacetime is
causal with closed $J^+(p)$, $J^-(q)$ for all $p$). A spacetime is
called globally hyperbolic if it is causal and the intersections
$J^+(p)\cap J^-(q)$ are compact for all $p, q \in M$ (for the last
two definitions, notice \cite{BeSa}).
 Globally hyperbolic spacetimes are the most relevant from both,
the geometric and physical viewpoints, and lie at the top of the
causal hierarchy .

\item A time or temporal function is called {\em Cauchy} if it is
onto on $\R$ and all its level hypersurfaces are Cauchy
hypersurfaces (i.e., topological hypersurfaces crossed exactly
once by any inextensible timelike curve). A classical theorem by
Geroch \cite{Ge} asserts the equivalence between: (i) to be
globally hyperbolic, (ii) to admit a Cauchy hypersurface, and
(iii) to admit a Cauchy time function. Moreover, the results in
\cite{BeSa1, BeSa2} also ensure the equivalence with: (iv) to
admit a (smooth) spacelike Cauchy hypersurface, and (v) to admit a
Cauchy temporal function ${\cal T}$. As a consequence, the full
spacetime admits a orthogonal Cauchy decomposition as in
(\ref{cosplit}).

\item Further
 properties have been achieved \cite{BeSa3}: any compact acausal
 spacelike submanifold with boundary can be extended to a (smooth) spacelike Cauchy
 hypersurface $\Sigma$, and any such $\Sigma$ can be regarded as
   a slice ${\cal T}=$constant for a suitable Cauchy orthogonal
 decomposition
 (\ref{cosplit}).
  Apart from the obvious
interest in the foundations of classical General Relativity, such
results have
 applications in fields such as the
wave equation or quantization, see for example \cite{BNP, Ru}.\eit
The following simple  results are useful for the discussions
below.

\bprop \label{p0} Let $(M,g)$ be a spacetime.

(1) If $\tau$ is a temporal function then there exists a conformal
metric $g^*= \Omega g$, $\Omega>0$, such that $\tau$ is steep.

(2) If ${\cal T}$ is a Cauchy temporal function and $\tau$ is a
temporal function then ${\cal T} + \tau$ is a Cauchy temporal
function. Moreover, ${\cal T} + \tau$ is steep if so is either
$\tau$ or ${\cal T}$.
 \eprop

\proof (1) As $\nabla^* \tau = \nabla \tau/\Omega$,  choose any
$\Omega \leq |\nabla \tau|^2$.


(2) ${\cal T} + \tau$ is temporal (and steep, if so is any of the
two functions) because of the reversed triangle inequality. In
order to check that its level hypersurfaces are Cauchy, consider
any future-directed timelike curve $\gamma: (a_-,a_+)\rightarrow
M$. It is enough to check that $\lim_{s\rightarrow \pm a} ({\cal
T} + \tau)(\gamma(s))= \pm \infty$. But this is obvious, because
$\lim_{s\rightarrow \pm a} {\cal T}(\gamma(s))= \pm \infty$ (as
${\cal T}$ is Cauchy) and $\tau(\gamma(s))$ is increasing.
 \cvd

\section{Characterization of isometrically embeddable Lorentzian
manifolds}\label{s3}

\bprop \label{t2.1} Let $(M,g)$ be a Lorentzian manifold. If there
exists a conformal immersion $i:M\rightarrow \LL^N$  then $(M,g)$
is a stably causal spacetime.

Moreover, if $i$ is a isometric immersion, then: (a) the natural
time coordinate $t=x^0$ of $\LL^N$ induces a steep temporal
function on $M$, and (b) the time-separation $d$ of $(M,g)$ is
finite-valued.

\eprop

\proof Notice that $x^0\circ i$ is trivially  smooth and also a
time function (as $x^0$ increases on $i\circ \gamma$, where
$\gamma$ is any future-directed causal curve in $M$), which proves
stable causality.

If $i$ is isometric, then $|\nabla (x^0\circ i)|\geq 1$ because,
at each $p\in M$, $\nabla (x^0\circ i)_p$ is the projection of
$\nabla x^0_{i(p)}$ onto the tangent space $di(T_pM)$, and its
orthogonal $di(T_pM)^\perp$ in $T_{i(p)}\LL^N$ is spacelike. This
 proves (a), for (b) notice that the finiteness of $d$ is an immediate
consequence of the
 finiteness of the time-separation $d_0$ on $\LL^N$ and the
straightforward inequality $d(p,q)\leq  d_0(i(p),i(q))$ for all
 $p,q\in M$.
 \cvd

\brema \label{r3a} {\em As a remarkable difference with the
Riemannian case,  Proposition \ref{t2.1} yields obstructions for
the existence of both, conformal  and isometric immersions in
$\LL^N$. In particular, non-stably causal spacetimes cannot be
conformally immersed, and further conditions on the
time-separation  are required for the existence of an isometric
immersion. In fact, it is easy to find even causally simple
spacetimes splitted as in (\ref{cosplit}) (with levels of ${\cal
T}$ non-Cauchy) which cannot be isometrically immersed in $\LL^N$,
as the following example shows. } \erema

\bexam \label{ex} {\em Let $M=\{(x,t)\in \R^2: x>0\}$, $
g=(dx^2-dt^2)/x^2$. This  is conformal to $\R^+\times \R \subset
\LL^2$ and, thus,  causally simple. 
 It is easy to check that $d(p,q)= \infty$ for
$p=(1,-2), q=(1,2)$ (any sequence of causal curves
$\{\gamma_m\}_m$ connecting $p$ and $q$ whose images contain
$\{(1/m, t): |t|<1\}$ will have diverging lengths). Thus, $(M,g)$
cannot be isometrically immersed in $\LL^N$.

 Recall that this example can be generalized, taking into account
that a stably causal spacetime is non-globally hyperbolic if and
only if it is conformal to a spacetime with a infinite-valued
time-separation (this holds for all strongly causal spacetimes,
see \cite[Th. 4.30]{BEE}). So, in the conformal class of any
non-globally hyperbolic spacetime, there are spacetimes non
isometrically inmersable in $\LL^N$.
 } \eexam
Nash's theorem will be essential for the proof of the following
result.

\bprop \label{t4} If a spacetime $(M,g)$ admits a steep temporal
function $\tau$ then it can be isometrically embedded in $\LL^N$
for some $N$. \eprop

\noi For the proof, recall first.

\blemm \label{l1} If a spacetime $(M,g)$ admits a  temporal
function $\tau$ then the metric $g$ admits a orthogonal
decomposition \be \label{barsplit} g= -\beta d\tau^2 + \bar g \ee
where $\beta=|\nabla \tau|^{-2}$ and $\bar g$ is a positive
semi-definite metric on $M$ with radical spanned by $\nabla \tau$.

In particular, if $\tau$ is steep then $\beta\leq 1$. \elemm

\proof The orthogonal decomposition (\ref{barsplit}) follows by
taking $\bar g$ as the trivial extension of $g|_{(\nabla
\tau)^\perp}$ to all $TM$. To determine the value of $\beta$,
recall that $d\tau(\nabla \tau) = g(\nabla \tau, \nabla \tau )=
-\beta \left(d\tau(\nabla \tau)\right)^2$ \cvd

{\em Proof of Proposition \ref{t4}.}
  Consider the orthogonal decomposition in Lemma \ref{l1}. Even
  though $M$  does not need to split as a product
  $\R \times
  S$ (in an open subset of $\LL^n$, the vector field $\nabla \tau$
may be incomplete and the topology of the level sets may change),
we can rewrite (\ref{barsplit}) as \be \label{tausplit} g= -\beta
d\tau^2 + g_\tau , \ee where each $g_{\tau_0}$ is Riemannian
metric on the slice $S_{\tau_0}=\tau^{-1}(\tau_0)$ varying
smoothly with $\tau_0$. Moreover, each $p\in M$ will be written as
$(\tau,x)$ where $x\in S_{\tau(p)}$.

Now, consider the auxiliary Riemannian metric
$$
 g_R:= (4-\beta) d\tau^2 + g_\tau .$$ By Nash's theorem, there
exists an isometric embedding $i_{\hbox{\small{nash}}}:(M,
g_R)\hookrightarrow \R^{N_0}$. Then, a simple computation shows
that the required isometric embedding $i: (M, g)\hookrightarrow
\LL^{N_0+1}$ is just:
$$
i(\tau,x)= \left( 2\tau, i_{\hbox{\small{nash}}}(\tau,x)\right).
$$
 \cvd

 \brema \label{r4} {\em (1) From the proof, it is clear that the
 hypotheses on steepness  can be weakened just by assuming
 that $\nabla \tau$ is lower bounded by some positive function
 $\epsilon(\tau)>0$.
 In fact, this is equivalent to require
 $\beta(\tau,x)\leq A(\tau)^2:=1/\epsilon(\tau)$, and the proof would work by taking
 $g_R:= (4A(\tau)^2-\beta) d\tau^2 + g_\tau $ and $
i(\tau,x)= \left( 2\int_0^\tau A(s)ds,
i_{\hbox{\small{nash}}}(\tau,x)\right)$. Nevertheless, no more
generality would be obtained in this case, because of the
following two different arguments: (a) it is easy to check that,
if this weaker condition holds, then a suitable composition $\hat
\tau = f\circ \tau$ for some increasing function $f$ on $\R$ would
be steep and temporal, and (b) the existence of a  steep temporal
function would be ensured by taking the isometric embedding $i:
M\hookrightarrow \LL^N$ and restricting the natural coordinate
$t=x^0$ as in Proposition \ref{t2.1}.

(2) Notice that Proposition \ref{t2.1} yields a necessary
condition for the existence of a isometric embedding and
Proposition \ref{t4} a sufficient one. Both together prove
trivially Theorem \ref{t1}, as well as Corollary \ref{c1}(2)
(notice also Proposition \ref{p0}(1)). Recall that, as a
difference with Nash's theorem, Proposition \ref{t4} does {\em
not} allow to prove that the spacetime is isometrically embedded
in an {\em arbitrarily small} open subset, which cannot be
expected now  (notice that $d(p,q) \leq d_0(i(p),i(q))\leq
|x^0(i(p))-x^0(i(q))|$). }
 \erema

\section{The Cauchy orthogonal b-decomposition of any globally hyperbolic
spacetime}\label{s4}

In order to obtain a steep Cauchy temporal function in a globally
hyperbolic spacetime, Proposition \ref{p0}(2) reduces the problem
to find a steep temporal function (not necessarily Cauchy), as the
existence of a Cauchy temporal function is ensured in
\cite{BeSa2}. Nevertheless, we will prove directly the existence
of  a  steep Cauchy temporal function ${\cal T}$, proving Theorem
\ref{t2} with independence of the results in \cite{BeSa2} (recall
Remark \ref{r}).

So, in what follows $(M,g)$ will be a globally hyperbolic
spacetime, and we will assume that $t$ is a Cauchy time function
as given by Geroch \cite{Ge}. The following notation will be also
used here. Regarding $t$,
$$
T^b_a= t^{-1}([a,b]), \quad \quad S_a= t^{-1}(a).$$
For any $p\in
M$, $j_p$ is the function
$$q\mapsto j_p(q) =\exp({-1/d(p,q)^2}).$$ For any $A,B\subset M$,
$$J(A,B):=J^+(A)\cap J^-(B)$$  in particular $J(p,S):=J^+(p)\cap
J^-(S)$ for $S$ any (Cauchy) hypersurface.

\subsection{Some  technical elements}\label{s4a}

In the next two propositions we will introduce a pair of technical
tools for the proof. But, first, consider the following
straightforward lemma, which will be invoked several times.

\blemm  \label{lemilla0} Let $\tau$ be a function  such that
$g(\nabla \tau, \nabla \tau) <0$ in some open subset $U$ and let
$K\subset U$ compact. For any function $f$ there exists a constant
$c$ such that $g(\nabla (f+c\tau), \nabla (f+c\tau)) <-1$ on $K$.

\elemm

\proof Notice that at each $x$ in the compact subset $K$ the
quadratic polynomial $g(\nabla (f(x)+c\tau(x)), \nabla
(f(x)+c\tau(x)))$ becomes smaller than -1
for some large $c$.
\cvd

 The following ``cone semi-time function'' will be useful from a technical viewpoint.

\bprop  \label{lemilla} Let $S$ be a Cauchy hypersurface, $p\in
J^-(S)$. For all neighborhood $V$ of $J(p,S)$ there exists a
smooth function $\tau\geq 0$ such that:

(i) $\hbox{{\rm Supp}} \, \tau\subset V$

(ii) $\tau > 1$ on $S\cap J^+(p)$.

(iii) $\nabla \tau$ is timelike and past-directed in $\hbox{{\rm
Int(Supp} } (\tau) \cap J^-(S))$.

(iv) $g(\nabla \tau, \nabla \tau) <-1$ on $J(p,S)$.

\eprop

{\em Proof.} Let $t$ be a Cauchy time function such
that\footnote{Along the proof, we will use this lemma only for
Cauchy hypersurfaces which are slices of a prescribed time
function. However, any Cauchy hypersurface can be written as such
a slice for some Cauchy time function. In fact, it is easy to
obtain a proof by taking into account that both, $I^+(S)$ and
$I^-(S)$ are globally hyperbolic and, thus, admit a Cauchy time
function --for details including the non-trivial case that $S$ is
smooth spacelike and $t$ is also required to be temporal, see
\cite{BeSa3}).}
 $S=S_a:=t^{-1}(a)$, and let $K\subset V$ be
a compact subset such that $J(p,S_a) \subset $ Int $(K)$.
Compactness guarantees the existence of  some $\delta >0$ such
that: for every $x\in K$ there exists a convex neighborhood
$U_x\subset V$ with $\partial^+U_x \subset J^+(S_{t(x)+2\delta})
$, where $\partial^+U_x:= \partial U_x \cap J^+(x)$. Now, choose
$a_0<a_1:= t(p) < \dots < a_n = a$ with $a_{i+1}-a_i<\delta /2$,
and construct $\tau$ by induction on $n$ as follows.

For $n=1$, cover $J(p,S)=\{p\}$ with a set type $I^+(x)\cap U_x$
with $x\in K\cap T_{a_0}^{a_1}$ and consider the corresponding
function $j_x$. For a suitable constant $c>0$, the product $c j_x$
satisfies both, (ii), (iii) and (iv). To obtain smoothability
preserving (i), consider the open covering $\{I^-(S_{a+\delta}),
I^+(S_{a+\delta/2})\}$ of $M$, and the first function $0\leq
\mu\leq 1$ of the associated partition of unity (Supp $\mu \subset
I^-(S_{a+\delta})$). The required function is just $\tau = c\mu
j_x$.

Now, assume by induction that the result follows for any chain
$a_0<\dots < a_{n-1}$. So, for any $k\leq n-1$, consider
$J(p,S_{a_k})$ and choose a compact set $\hat K \subset $ Int $K$
with $J(p,S)\subset $ Int $\hat K$. Then, there exists  a function
$\hat \tau$ which satisfies condition (i) above for $V=$ Int $\hat
K \cap I^-(S_{a_{k+1}})$ and conditions (ii), (iii), (iv) for
$S=S_{a_k}$. Now, cover $\hat K \cap T^{a_{k+1}}_{a_k}$ with a
finite number of sets type $I^+(x^i)\cap U_{x^i}$ with $x^i\in K
\cap T^{a_{k+1}}_{a_{k-1}}$, and consider the corresponding
functions $j_{x^i}$.

For a suitable constant $c>0$, the sum $\hat \tau + c \sum_i
j_{x^i}$ satisfies  (iii) for $S=S_{a_{k+1}}$. This is obvious in
$J^-(S_{a_k})$ (for any $c>0$), because of the convexity of
timelike cones and the reversed triangle inequality. To realize
that this can be also obtained in $T^{a_{k+1}}_{a_k}$, where $
\nabla \tau$ may be non-timelike, notice that the support of
$\nabla \hat\tau |_{T^{a_{k+1}}_{a_k}}$ is compact, and it is
included in the interior of the support of $\sum_i j_{x^i}$, where
the gradient of the sum is timelike; so, use Lemma \ref{lemilla0}.
 As $J(p,S_{a_{k+1}})$ is compact, conditions  (ii), (iv) can be
trivially obtained by choosing, if necessary, a bigger $c$.

Finally, smoothability (and (i)), can be obtained again by using
the open covering $\{I^-(S_{a_{k+1}+\delta})$, $
I^+(S_{a_{k+1}+\delta/2})\}$ of $M$, and the corresponding first
function $\mu$ of the associated partition of  unity, i.e. $\tau =
\mu (\hat \tau + c \sum_i j_{x^i})$. \cvd

In order to extend  locally defined time functions to a global
time one, one cannot use a partition of  unity (as stressed in the
previous proof, because $\nabla \tau$ is not always timelike when
$\mu$ is non-constant). Instead,  local time functions must be
added directly and, then, coverings as the following  will prove
useful.

\bdefi Let $S$ be a Cauchy hypersurface. A {\em fat cone covering}
of $S$ is a sequence of  pairs of  points $p'_i\ll p_i , i\in \N$
such that both\footnote{Strictly, we will need only the local
finiteness of ${\cal C}'$.}, ${\cal C}'=\{I^+(p'_i): i\in \N\}$
and ${\cal C}=\{I^+(p_i): i\in \N\}$ yield a locally finite
covering of $S$.
 \edefi


 \bprop \label{fatcone} Any Cauchy hypersurface $S$ admits a fat cone covering $p'_i \ll p_i, i\in \N$.

 Moreover, both ${\cal C}$ and ${\cal C}'$ yield also a finite
 subcovering of $J^+(S)$.
\eprop \proof  Let $\{K_j\}_j$ be a sequence of compact subsets of
$S$ satisfying $K_j\subset $ Int $K_{j+1}$, $S=\cup_jK_j$. Each
$K_j\backslash$ Int $K_{j-1}$ can be covered by a finite number of
sets type $I^+(p_{jk}), k=1\dots k_j$ such that $I^+(p_{jk})\cap S
\subset K_{j+1}\backslash K_{j-2}$. Moreover, by continuity of the set-valued function $I^+$, this last inclusion
is fulfilled if each $p_{jk}$ is replaced by some close
$p'_{jk}\ll p_{jk}$, and  the required pairs $p'_i (=p'_{jk})$,
$p_i (=p_{jk})$, are obtained.

For the last assertion, take $q\in J^+(S)$ and any compact
neighborhood $W \ni q$. As $J^-(W)\cap S$ is compact, it is
intersected only by finitely many elements of ${\cal C}, {\cal C}'$,
and the result follows. \cvd

\subsection{Construction of the b-decomposition}\label{s4b}

\begin{definition}
\label{fls} Let $\q , p \in T_{a-1}^a$, $ \q\ll p$. A {\em steep
forward cone  function (SFC)} for $(a, p' , p) $ is a smooth
function $ h_{a,\q,p}^+ : M \rightarrow [ 0 , \infty )$ which
satisfies the following:

\begin{enumerate}
\item{$Supp (h_{a,\q,p} ^+)  \subset J(\q, S_{a+2})$,}
 \item{$h_{a,\q,p}^+ > 1$ on $
S_{a+1}\cap J^+ (p)$},
 \item{If $x
\in J^- (S_{a+1}) $ and $h_{a,\q,p}^+ (x) \neq 0 $ then $ \nabla
h_{a,\q,p} ^+ (x)$ is timelike and past-directed, and} \item{ $g(
\nabla h_{a,\q,p}^+ , \nabla h_{a,\q,p} ^+ ) < - 1 $ on
$J(p,S_{a+1})$.}
\end{enumerate}

\end{definition}

Now,  Proposition \ref{lemilla} applied to $S=S_{a+1},
V=I^-(S_{a+2})\cap I^+(\q)$ yields directly:

\bprop For all $(a,\q,p)$ there exists a SFC. \eprop

The existence of a fat cone covering (Proposition \ref{fatcone})
allows to find a function $h^a_+$ which in some sense globalizes
the properties of a SFC.

\blemm \label{l07} Choose  $a\in \R$ and take any fat cone
covering $ \{ p'_i\ll p_i \vert i\in \N \}$ for $S=S_a$. For every positive
sequence $\{ c_i \geq 1 \vert i \in \N \}$, the non-negative
function $h_{a}^+:=(|a|+1) \sum_i c_i h_{a,\q_i,p_i} ^+$
satisfies:

\begin{enumerate}
\item{$Supp (h_{a} ^+)  \subset J(S_{a-1},S_{a+2})$,}
 \item{$h_{a}^+ > |a|+1$ on\footnote{This condition is imposed in order to ensure that the finally obtained  temporal function is Cauchy. It could be dropped
 if one  looks only for a temporal function and, then, uses Proposition \ref{p0}(2).}
 $S_{a+1}$},
 \item{If $x
\in J^- (S_{a+1}) $ and $h_{a}^+ (x) \neq 0 $ then $ \nabla h_{a}
^+ (x)$ is timelike and past-directed, and} \item{ $g( \nabla
h_{a}^+ , \nabla h_{a} ^+ ) < - 1 $ on $J(S_a, S_{a+1})$.}
\end{enumerate}
\elemm \proof Obvious. \cvd

The gradient of $h_a^+$ will be spacelike  at some subset of
$J(S_{a+1}, S_{a+2})$. So, in order to carry out the inductive
process which proves Theorem \ref{t2}, a strengthening of Lemma
\ref{l07} will be needed.

\blemm \label{l07b} Let $h^+_a\geq 0$ as in Lemma \ref{l07}. Then
there exists a function $h^+_{a+1}$ which satisfies all the
properties
 corresponding to Lemma \ref{l07} and additionally:

\be \label{eee} g(\nabla (h^+_{a}+  h^+_{a+1}), \nabla (h^+_{a}+
h^+_{a+1}))< -1 \quad \quad \hbox{on} \; J(S_{a+1}, S_{a+2}) \ee
 (so, this inequality holds automatically on all $ J(S_{a}, S_{a+2})$).
\elemm

\proof Take a fat cone covering $ \{ p'_i\ll p_i \vert i\in \N \}$
for $S=S_{a+1}$. Now, for each $p_i$ consider a constant $c_i\geq
1$ such that $c_i h_{a+1,\q_i,p_i} ^+ + h^+_a$ satisfies
inequality (\ref{eee})  on $J(p_i, S_{a+2})$ (see Lemma
\ref{lemilla0}). The required function is then $h_{a+1}^+=(|a|+2)
\sum_i c_i  h_{a+1,\q_i,p_i} ^+$.
 \cvd

Now, we have the elements to complete our main proof.

\noi {\em Proof of Theorem \ref{t2}.} Consider the function
$h^+_a$ provided by Lemma \ref{l07} for $a=0$, and apply
inductively Lemma \ref{l07b} for $a=n\in \N$. Then, we obtain a
function ${\cal T}^+ = \sum_{n=0}^\infty h_n^+\geq 0$ with nowhere
spacelike gradient, which is a steep temporal function on
$J^+(S_0)$ with support in $J^+(S_{-1})$. Analogously, one can
obtain a function ${\cal T}^-\geq 0$ which is a steep temporal
function with the reversed time orientation, on $J^-(S_0)$. So,
${\cal T}={\cal T}^+-{\cal T}^-$ is clearly a steep temporal
function on all $M$.

Moreover, the levels hypersurfaces of ${\cal T}$ are Cauchy. In
fact, consider  any future-directed causal curve  $\gamma$, and
reparametrized it with the Cauchy time function $t$. Then,
$$
\lim_{t\rightarrow \infty} {\cal T}(\gamma(t)) \left( = \lim_{n\in
\N} {\cal T}^+(\gamma(n+1))\geq \lim_{n\in \N}
h^+_n(\gamma(n+1))\right) = \infty, \quad \lim_{t\rightarrow
-\infty} {\cal T}(\gamma(t))=-\infty , $$ and $\gamma$ crosses all
the levels of ${\cal T}$, as required. \cvd

\section{Appendix}

Clarke \cite{Cl} developed the following method in order to embed
isometrically any manifold $M$ endowed with a semi-Riemannian (or
even degenerate)  metric $g$ in some semi-Euclidean space
$\R^N_s$. First, he proved that, for some $p\geq 0$, there exists
a function $f: M\rightarrow \R^p_p$ such that the  (possibly
degenerate) pull-back metric $g(f)$ on $M$ induced from $f$
satisfies $g_R= g-g(f)>0$. So, the results for positive definite
metrics are applicable to $(M,g_R)$, and one can construct a
Riemannian isometric embedding $f_R: M\rightarrow \R^{N_0}$ ($f_R$
can be constructed from Nash result, even though Clarke develops a
technique
to reduce the Nash value for $N_0$). Then, the required embedding
$i: M \rightarrow \R^N_p$ is obtained as a product $i(x)= (f(x),
f_R(x))$ for $N=p+N_0$.

\sm \sm

\noindent In Lorentzian signature, Clarke's optimal value for $p$
is 2. Nevertheless, he claims that, if $(M,g)$ is a globally
hyperbolic spacetime, then one can take $p=1$ \cite[Lemma 8]{Cl}.
Our purpose in this Appendix is to analyze this  question and
show:

(A) the required condition $g-g(f)>0$ on $f$   is essentially
equivalent to be a steep temporal function, and

(B) the success of the construction of $f$ in \cite[Lemma 8]{Cl}
depends on a  problem of smoothability, which may have interest in
its own right.

\sm \sm

\noindent In order to make  these points clear, we will
particularize the proof of \cite[Lemma 8]{Cl} to a very simple
case, and will follow most of the notation there. As a previous
remark, Clarke assumed that the existence of a temporal function
$\tau$ had already been  proved, as this question (one of the
prominent folk problems of smoothability) seemed true then. At any
case, we can assume now even that $\tau$ is Cauchy temporal. Then,
consider a globally hyperbolic spacetime which can be written as
$$(\R^2,g) \quad \quad g=-V^2 d\tau^2 + M^2 dy^2 , $$
where $(\tau, y)$ are the natural coordinates of $\R^2$ and $V, M$
are two positive functions on $\R^2$. Easily, a function $f: \R^2
\rightarrow \R^1_1 (=\LL^1)$ satisfies $g-g(f)>0$ if and only if:
\be \label{cl1} -V^2 (\partial_y f)^2 + M^2 (\partial_\tau f)^2 >
V^2M^2, \ee and this is trivially equivalent to $g(\nabla f,
\nabla f) < -1$. This proves (A) in our particular example and,
taking into account Remark \ref{r4}(1), it seems general.

Now, consider any smooth function $\sigma\geq 0$ on $\R^2$
invariant through the flow of $\nabla \tau$ such that
$\sigma^{-1}([0,s])$ is compact for all $s$, and let
$Y=\sigma^{-1}([0,1])$; in our simplified example, we can put
$\sigma(\tau, y)=|y|^2$. Outside $Y$ the two lightlike vector
fields,
$$ A_\pm = M
\partial_\tau \pm V \partial_\sigma , $$
are well defined, and equation (\ref{cl1}) can be also rewritten
as \be \label{c2} (A_+f)(A_-f)> V^2 M^2. \ee So, the crux  is to
construct a function $f$ which  satisfies (\ref{c2}) outside $Y$,
among other conditions. Clarke's proposal is the following. Let
$$
H^\pm (t, s) = J^\pm(\tau^{-1}(0)) \cap J^\mp(\tau^{-1}(t)\cap
\sigma^{-1}([0,s])). $$ After choosing a certain volume element
$\omega$, the function $f$ is defined  as: \be\label{c3} f(x)=
\int_{H^+ (\tau(x), \sigma(x))} \omega \ee  whenever
$\tau(x)>\epsilon>0$ and outside a neighborhood of\footnote{For
$\tau(x)<-\epsilon<0$, the function $f$ is negative and defined
dually in terms of $H^-$, for $\tau(x)=0$, $f$ is $0$, and a more
technical definition is given for $f$ on a neighborhood of $Y \cup
\tau^{-1}(0)$. However, this is not relevant for our discussion. }
$Y$. Notice that $A_\pm$ are future directed, and $A_+$ points
outwards the region $\sigma^{-1}([0,\sigma(x)])$ at each $x\in
M\backslash Y$. So, if
 $f$ is $C^1$, then one would expect $A_+(f)>A_-(f)>0$. Moreover, Clarke claims that (\ref{c2}) can be
 also achieved by choosing
  $\omega$ large enough (and   eventually, a redefinition of $\tau$).

 At what extent can one assume that $f$ is $C^1$ (or, at least,
 that it can be smoothed to a function which satisfies the
 required conditions)? For each measurable subset $Z$ of the spacetime manifold,
 consider its $\omega$-measure
 $\mu(Z)=\int_Z \omega$. In any causally continuous spacetime $M$ it
 is known that the functions $x\mapsto \mu (J^\pm(x))$ are
 continuous, if $\mu(M)<\infty$. Moreover, if $M$ is globally
 hyperbolic and $S$ is any topological Cauchy hypersurface,
then $I^+(S)$ is a globally hyperbolic spacetime in its own right,
and the function $x\mapsto \mu (J(S, x)), x\in I^+(S)$, becomes
continuous, even if we drop the assumption about the finiteness of
$\mu$. Nevertheless, neither the functions $\mu (J^\pm(x))$ nor
$\mu (J(S, x))$ are smooth in general (see figure). In Clarke's
case, the fact that $S=\tau^{-1}(0)$ is not only smooth but
spacelike, may help to stablish smoothness. However, recall that
the definition of $f$ also uses the function $\sigma$. Such a
$\sigma$ can be defined by taking some auxiliary complete
Riemannian metric on $S$, and smoothing along the cut locus the
squared distance function  to a fixed point $y_0\in S$. The
behavior of $f$ at the points $x\in M$ such that the boundary of
$S\cap J^\mp(\tau^{-1}(x)\cap \sigma^{-1}([0,\sigma(x)]))$
intersects the cut locus may complicate the situation.

Summing up, even assuming --as a necessary element of Clarke's
proof-- the existence of a  temporal function, which was proved
in \cite{BeSa2} and is re-proved in a shorter form here, the
smoothability of $f$ remains as a non-trivial  problem. The
solution of this question not only would complete Clarke's proof
but also may have interest in its own right.

\begin{figure}[ht]
\centering
\includegraphics[width=12cm]{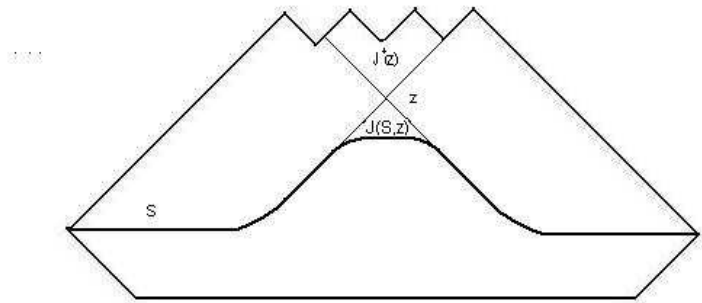}
\caption{The depicted open subset of $\LL^2$ is globally
hyperbolic, and $S$ a smooth Cauchy hypersurface. Functions
$J^+(x)$ and $J(S,x)$ are not smooth at $z \in I^+(S)$. }
\label{fig1}
\end{figure}

\section*{Acknowledgements}

 The  comments by Prof. L. Andersson, who emphasized reference
 \cite{Cl}, and the careful reading by the referee,
 are deeply acknowledged.

MS is partially supported by  Grants P06-FQM-01951 (J.
Andaluc\'{\i}a) and  MTM2007-60731 (MEC with FEDER funds). OM was
partially funded by CoNaCyT projects 49093 and 82471.


\end{document}